\documentclass[12pt]{amsart}
\usepackage{amssymb,amscd,amsfonts,amsbsy,amsmath,amsthm,amsrefs, mathtools}
\usepackage{graphicx}
\usepackage{xcolor}
\usepackage{caption}
\usepackage{subcaption}
\usepackage{pdfsync}	
\usepackage{dsfont}
\usepackage{keyval} 
\usepackage{float}				
\usepackage{xcolor}
\usepackage[unicode=true,psdextra]{hyperref}
\hypersetup{
	colorlinks   = true,    
	urlcolor     = blue,    
	linkcolor    = blue,    
	citecolor    = red      
}
\usepackage[scr=boondox]{mathalpha}
\usepackage{upgreek}
\usepackage{enumitem}

\newcommand{\dive}{\operatorname{div}}

\newcommand{\dy}{\,\text{d}y}
\renewcommand{\d}{\, \text{d}}

\definecolor{bostonuniversityred}{rgb}{0.8, 0.0, 0.0}

\usepackage{etoolbox}
\usepackage[letterpaper, left=1.9cm, right=1.9cm, top=2.5cm, bottom=3cm, footskip=1.5cm]{geometry}
\patchcmd{\section}{\scshape}{\bfseries\Large}{}{}
\makeatletter
\renewcommand{\@secnumfont}{\bfseries\Large}
\def\l@subsection{\@tocline{2}{0pt}{4pc}{5pc}{}}
\makeatother

\setlist[description]{leftmargin=.4cm,labelindent=0cm}

\newcommand{\LL}{\mathbb{L}}
\newcommand{\R}{{\mathbb{R}}}
\renewcommand{\P}{{\mathbb{P}}}
\newcommand{\N}{{\mathbb{N}}}
\newcommand{\CC}{{\mathbb{C}}}

\newcommand{\BBc}{\mathcal{B}}

\newcommand{\esssup}[1]{\underset{#1}{\text{\hspace{1pt}ess\hspace{2pt}sup }}} 

\newtheorem{theorem}{Theorem}[section]
\newtheorem{proposition}[theorem]{Proposition}
\newtheorem{lemma}[theorem]{Lemma}

\theoremstyle{definition}

\newtheorem{remark}[theorem]{Remark}
\numberwithin{equation}{section}

\usepackage{framed,color}
\definecolor{shadecolor}{rgb}{0.9,0.9,0.9}

\title[Geometric Structures in 2D MHD Flows]{Spontaneous Geometric Structures in 2D MHD Flows}
\author[M. Nasi]{Matheus Nasi} 
\address[M. Nasi]{Institut Camille Jordan, 
	Université de Lyon, Université Claude Bernard Lyon~1, 
	69622 Villeurbanne Cedex. France and Departamento de Matem\'atica. 
	Universidade Federal de Pernambuco, Recife, PE 50740-560. Brazil}
\email{matheus.nasi@ufpe.br}
\email{carneiro@math.univ-lyon1.fr}
\author[C. F. Perusato]{Cilon F. Perusato}
\address[C. F. Perusato]{Departamento de Matem\'atica. 
	Universidade Federal de Pernambuco, Recife, PE 50740-560. Brazil 
	}
\email{perusato@math.univ-lyon1.fr}
\email{cilon.perusato@ufpe.br}

\begin{document}
	
	\date{}
	
	\maketitle

	\pretolerance=10000
	\begin{abstract}
		Under suitable spatial decay assumptions on the initial data, we establish precise asymptotic expansions of the velocity and magnetic fields for two-dimensional incompressible MHD flows at each positive time. We identify in the far-field profile the hexagonal pattern previously discovered for the 2D Navier–Stokes equations by L. Brandolese \cite{Bra22}, characterized by six equally spaced directions of faster decay for each component of  $u(x,t)-u_0(x)$. Surprisingly, the corresponding contribution to the magnetic field vanishes identically, owing to an exact cancellation associated with the antisymmetric structure of the magnetic nonlinearity. Thus, at order $|x|^{-3}$, the hexagonal structure appears only in the velocity field.
	\end{abstract}
\section{Introduction}
Two-dimensional Navier–Stokes flows can spontaneously develop hexagonal structures in the far field, even when no symmetry is imposed on the initial data \cite{Bra22}. Does this geometry persist in magnetohydrodynamic flows, and does it appear in both the velocity and magnetic fields? In this paper, we show that the two fields behave differently at order $|x|^{-3}$. 

The magnetohydrodynamics equations can be written in the following way
\begin{equation}
	\label{eq:MHD}
	\begin{cases}
		\partial_t u + (u \cdot \nabla)u - S(b \cdot \nabla)b + \nabla\left(p + \frac{S}{2}|b|^2\right) = \frac{1}{R_e}\Delta u \\
		\partial_t b + (u \cdot \nabla)b - (b \cdot \nabla)u = \frac{1}{R_m}\Delta b \\
		\nabla\cdot u = \nabla\cdot b = 0 \\
		u(\cdot, 0) = u_0, \, b(\cdot,0) = b_0.
	\end{cases}.
\end{equation}
Here, the unknowns are the velocity field $u$ of the fluid, the pressure $p$ and the magnetic field $b$, all defined in $\R^2$. The positive constants $R_e$ and $R_m$ are, respectively, the Reynolds number and the magnetic Reynolds number. Moreover, $S = M^2/(R_e R_m)$, where $M$ is the Hartmann number. After rescaling $u$ and $b$, we can assume $S = R_e = 1$, and, to simplify the calculations, we shall assume $R_m = 1$.

Applying the Leray-Hopf projector $\P$ and then the Duhamel formula, we obtain the integral equations
\begin{equation}\tag{MHDI}
	\label{eq:MHDintegral}
	\begin{cases}
		u(t) = e^{t\Delta}u_0 - \int_0^t e^{(t-s)\Delta}\P\nabla\cdot ( u \otimes u - b \otimes b)(s) \d s \\
		b(t) = e^{t\Delta}b_0 - \int_0^t e^{(t-s)\Delta}\P\nabla\cdot (u \otimes b - b \otimes u)(s) \d s \\
		\nabla\cdot u_0 = \nabla\cdot b_0 = 0
	\end{cases}.
\end{equation}
If $(u_0,b_0) \in L_{\sigma}^2 \times L_{\sigma}^2$, then this system is known to possess a unique global Leray weak solution (see \cite{Leray34,SM83}), that is, a weak solution to (\ref{eq:MHD}), such that
\[
(u,b) \in \left(L^{\infty}((0,\infty); L^2_{\sigma}(\R^2))\cap L^2((0,\infty); \dot{H}^1(\R^2))\right)^{2 \times 2}, 
\]

and satisfying the energy inequalities, for all $t > 0$,
\begin{equation}
	\label{ineq:MHDEnergy}
	\|u(t)\|_2^2 +  \|b(t)\|_2^2 + 2 \int_0^t \|\nabla u(s)\|_2^2 + \|\nabla b(s)\|_2^2 \d s \leq \|u_0\|_2^2 + \|b_0\|_2^2.
\end{equation}

Roughly, the results of \cite{Bra22} describe the behavior of two-dimensional {\it Navier–Stokes flows} in the far field as follows. Let $v(\cdot,t)$ denote any component of the velocity field $u(\cdot,t)$. Then, for each fixed $t>0$ and $|x|\gg 1$,  
\begin{enumerate}
	\item[(i)] If $u_0$ is {\it spatially well localized} in $\R^2$, the level lines of $|v(\cdot,t)|$ tend to have a hexagonal structure. 
	\item[(ii)] If $u_0$ is {\it spatially poorly localized} in $\R^2$, the level lines of $|v(\cdot,t) - v_0|$ tend to have a hexagonal structure. 
\end{enumerate}   
The proof of such a statement relies on a refined asymptotic formula for the flow in the far-field, namely

\begin{equation*}
\label{eq:NSAsym}
	u(x,t) = e^{t\Delta}u_0(x) - \frac{L_{NS}(t)}{|x|^3}\begin{pmatrix}
		\cos(3\theta + \alpha(t)) \\ \sin(3\theta + \alpha(t))
	\end{pmatrix}
	+ o_t(|x|^{-3}) \qquad \text{as } |x| \to \infty,
\end{equation*}
where $x=|x|(\cos\theta,\sin\theta)$ and $\alpha(t)$ is a phase. Here, \(o_t(|x|^{-3})\) denotes a remainder that is \(o(|x|^{-3})\) as \(|x|\to\infty\) for each fixed \(t>0\).
The amplitude $L_{NS}(t)$ is given by
\[
L_{NS}(t)
:=\lim_{|x|\to\infty}|x|^3|u(x,t)-u_0(x)|
=\frac{1}{\pi}
\left[
\left(\int_0^t\int_{\mathbb R^2}(u_1^2-u_2^2)\,dy\,ds\right)^2
+
\left(\int_0^t\int_{\mathbb R^2}2u_1u_2\,dy\,ds\right)^2
\right]^{1/2}.
\]
The hexagonal structures described above occur at times when $L_{NS}(t)>0$.

The MHD system, however, contains two coupled quadratic tensors with fundamentally different algebraic structures. Indeed, the velocity equation involves the symmetric stress  $u \otimes u - b \otimes b$, whereas the magnetic equation involves the antisymmetric tensor $u\otimes b-b\otimes u$. This distinction produces the main new phenomenon of the present work. In the velocity equation, the leading nonlinear contribution is generally nonzero and gives rise to the same third-order angular structure as in the planar Navier–Stokes case, with the kinetic stress replaced by the difference between the kinetic and magnetic stresses. In contrast, in the magnetic equation the corresponding term of order $|x|^{-3}$ vanishes identically.

The energy estimate \eqref{ineq:MHDEnergy} ensures that the following integrals are finite for every $t>0$. We can therefore associate with the solution $(u,b)$ the complex-valued curve

\begin{equation*}
	{z}(t) := \left( \int_0^t \!\! \int (u_1^2 - u_2^2 - b_1^2 + b_2^2) \ dy \hspace{1pt} ds \right) + i \left( \int_0^t \!\! \int 2(u_1u_2 - b_1b_2) \ dy \hspace{1pt} ds \right) ~ \in \CC.
\end{equation*}
The use of integrals of such form is a generalization to the MHD system inspired by earlier works done by Schonbek \cite{Sch91}, Dobrokhotov and Shafarevich \cite{DobSha} and Miyakawa and Schonbek \cite{MyiSch}.

We can now state our main result. 
\begin{theorem}
	\label{thm:MHDthm}
	Let $u_0, b_0 \in L_{\sigma}^2(\R^2) \cap L^p(\R^2)$ for some $4 < p <\infty$ be such that
	\begin{equation}
		\label{eq:InitialDecay}
		\begin{cases}
			u_0(x) = o(|x|^{-1}(\log|x|)^{-1/2}) \\
			\nabla u_0(x) = o(|x|^{-2}(\log|x|)^{-1/2}) \\
			\nabla^2 u_0(x) = o(|x|^{-3})
		\end{cases}
		\text{ and } \quad
		\begin{cases}
			b_0(x) = o(|x|^{-1}(\log|x|)^{-1/2}) \\
			\nabla b_0(x) = o(|x|^{-2}(\log|x|)^{-1/2}) \\
			\nabla^2 b_0(x) = o(|x|^{-3})
		\end{cases}
	\end{equation}
	as $|x| \to \infty$. Let $(u,b)$ be the unique global Leray solution starting from $(u_0,b_0)$. Then, for all $t > 0$, there exists a phase $\phi(t)$ such that
	\begin{equation}
	\label{eq:uProfile}
		u(x,t) = e^{t\Delta}u_0(x) + \frac{{L}(t)}{|x|^3}
		\begin{pmatrix}
			\cos (3\theta + \phi(t)) \\ \sin(3\theta + \phi(t))
		\end{pmatrix}
		+ o_t(|x|^{-3})
	\end{equation}
	and
	\begin{equation*}
		b(x,t) = e^{t\Delta}b_0(x) + o_t(|x|^{-3})
	\end{equation*}
	as $|x| \to \infty$. Here, ${L}(t)$ is given by
	\begin{equation*}
		{L}(t) := \frac{1}{\pi}|{z}(t)|.
	\end{equation*}
\end{theorem} 
When $L(t)>0$, each component of $u(\cdot,t)-u_0$ has exactly six equally spaced directions of faster spatial decay, whereas
\[
|u(x,t)-u_0(x)|=\frac{L(t)}{|x|^3}+o_t(|x|^{-3}).
\]
Thus, the components exhibit a hexagonal pattern, while the magnitude has a radial leading profile. The magnetic difference $b(\cdot,t)-b_0$ has no contribution at this order. If, in addition, $u_0(x),b_0(x)=o(|x|^{-3})$, these conclusions apply directly to the velocity and magnetic fields.

We note that the decay assumptions (\ref{eq:InitialDecay}) are natural for finite energy flows. Indeed, the first of (\ref{eq:InitialDecay}) is just a pointwise analogue of the condition $u_0,b_0 \in L^2(\R^2)$, and often $\nabla u_0$ decays one decay rate faster than $u_0$, which is enough for the second condition in (\ref{eq:InitialDecay}).

Discrete rotational symmetries also arise  for inviscid flows. For example, in \cite{CCJavi}, the authors constructed smooth, compactly supported vorticities for the two-dimensional Euler equations that rotate uniformly and possess $m$-fold symmetry, including the case $m=6$. Their construction concerns exact rotating solutions with symmetric vorticity level sets. The hexagonal structure considered here instead describes the leading far-field correction to the velocity and requires no rotational symmetry of the initial data. In the MHD setting, this asymptotic geometry is further distinguished by the cancellation of the corresponding contribution to the magnetic field.

The proof will follow these steps: based on the works of Kato \cite{Kato84}, we show the existence of a solution in a certain Kato-like space $(Y_{p,T_0})^2$, for $T_0>0$ small enough, via a fixed-point argument. Then, by a variation of \cite{Dan05}, we show that such a solution agrees with the unique \textit{Leray}-type global solution $(u,b)$ and, by adapting a bootstrapping argument in \cite{Bra22,Vig05}, we show that the spatial decay of both fields and their gradients is preserved. Once we have those estimates, we can use a decomposition of the kernel of $e^{t\Delta}\P\textmd{div}$ developed in \cite{Bra22} to obtain our desired profiles.

The paper is organized as follows. Section \ref{section2} collects preliminary estimates, and Section \ref{section3} contains the proof of the main theorem. The discussion is supplemented by an appendix, where we obtain lower bounds for the radius beyond which the hexagonal structures emerge.

\section{Preliminaries}\label{section2}
In this section, we present some initial results regarding the kernel of the bilinear operator associated to (\ref{eq:MHD}) that will be useful throughout the paper, especially in showing the existence of solutions to (\ref{eq:MHDintegral}). We also show a time-dependent pointwise control of $|e^{t^\Delta}(u_0,b_0) - (u_0, b_0)|$.
\subsection*{Decomposition of the kernel $e^{t\Delta}\P\textmd{div}$ in $\R^d$, $d \geq 2$}
Consider $F(x,t) = (F_{j;h,k})(x,t)$ the kernel of the operator $e^{t\Delta}\P\text{div}$. $F$ can be defined through its Fourier symbol
\begin{equation}
\label{def:OseenKernel}
	\hat{F}_{j,h,k}(\xi,t) = e^{-t|\xi|^2}\left( i\xi_h\delta_{j,k} + \frac{i\xi_j i\xi_h i\xi_k}{|\xi|^2}\right), \qquad j,h,k = 1, \ldots, d, \quad \xi \in \R^d, \quad t > 0,
\end{equation}
where $\delta_{j,k}$ is the Kronecker symbol. We can write the integral formulation of the MHD system in $\R^d$, $d \geq 2$, as

\begin{equation*}
	\begin{cases}
		u_j(x,t) = e^{t\Delta}u_{0,j}(x) - \int_0^t \int \sum_{h,k=1}^d F_{j;h,k}(x-y,t-s)(u_hu_k - b_hb_k)(y,s) \d y \d s \\
		b_j(x,t) = e^{t\Delta}b_{0,j}(x) - \int_0^t \int \sum_{h,k=1}^d F_{j;h,k}(x-y,t-s)(u_hb_k - b_hu_k)(y,s) \d y \d s
	\end{cases},
\end{equation*}

for $j = 1, \ldots, d$, together with the divergence-free conditions $\nabla \cdot u_0 = \nabla \cdot b_0 = 0$. Computing the inverse Fourier transform of (\ref{def:OseenKernel}) yields the decomposition $F = F^{(1)} + F^{(2)}$, with
\begin{equation*}
	F^{(1)}_{j;h,k}(x,t) = \delta_{j,k}\partial_hg_t \qquad F^{(2)}_{j,h,k}(x,t) = \int_t^{\infty} \partial_j \partial_h \partial_k g_s (x) \d s.
\end{equation*}

From this decomposition, we get that the kernel admits the pointwise inequality
\begin{equation*}
	|F(x,t)| \leq C_{\eta}|x|^{-\eta}t^{\frac{\eta - d - 1}{2}}, \qquad 0 \leq \eta \leq d + 1,
\end{equation*} 
for some $C_{\eta} > 0$. Moreover, we also have estimates on its $L^p$ norms
\[
	\|F(t)\|_p \leq C t^{-\frac{d+1}{2} + \frac{d}{2p}}, \qquad 1 \leq p \leq \infty.
\]

Another decomposition of the kernel $F$, due to \cite{Bra08}, is given by 
\begin{equation}
	\label{def:KernelDecomPBrandolese}
	F(x,t) = \mathfrak{F}(x) + |x|^{-d-1}\Psi\left(x/\sqrt{t}\right),
\end{equation}
where $\mathfrak{F}$ is a homogeneous tensor of degree $-d-1$, whose components are given by
\begin{equation*}
	\mathfrak{F}_{j;h,k}(x) = \partial_j \partial_h \partial_k E_d,
\end{equation*}
where $E_d$ is the fundamental solution of the Laplacian in $\R^d$. Moreover, $\Psi = (\Psi_{j;h,k})$ is smooth outside the origin and such that, for all $\alpha \in \N^d$, and $x \neq 0$, there exist $C, c > 0$ such that
\begin{equation*}
	|\partial^{\alpha}\Psi(x)| \leq Ce^{-c|x|^2}.
\end{equation*}
Such decomposition of $F$ proves to be very important in the study of far-field asymptotics of the bilinear term in (\ref{eq:MHDintegral}), due to (\cite{Bra22}, Lemma 5.1), since, by defining, for $t  > 0$,
\begin{equation*}
	\LL(w)(x,t) := \int_0^t \int F(x-y,t-s)w(y,s) \d y \d s,
\end{equation*}
we get the following decomposition
\begin{lemma}
	\label{lemma:5.1}
	Let $0 < T < \infty$, $w \in (L^1(\R^d \times (0,T)))^{d \times d}$ and $0 \leq \kappa < 1$. Assume that 
	\begin{equation}
	\label{eq:L5.1}
		\esssup{t \in (0,T)} t^{\kappa} |\nabla w(x,t)| = o(|x|^{-d-1} (\log |x|)^{-1}), \qquad \text{as } |x| \to \infty.
	\end{equation}
	Then, for $|x| \neq 0$ and $t \in (0,T)$,
	\begin{equation}\label{conclusion:L5.1}
		\LL(w)(x,t) = \mathfrak{F}(x) : \int_0^t \int w(y,s) \ dy \hspace{1pt} ds + |x|^{-d-1}\epsilon(x,t),
	\end{equation}
	where
	\begin{equation*}
		|\epsilon(x,t)| \to 0 \quad \text{as } |x| \to \infty.
	\end{equation*}
\end{lemma}
The colon symbol stands for a summation on the last two subscripts of $\mathfrak{F}_{j;h,k}$.
\begin{remark}\label{needforappendix}
	If we assume that $\textmd{ess\,sup}_{t\in(0,T)} t^{\kappa} |w(x,t)| = o(|x|^{-d-1})$ as $|x| \to \infty$ holds instead of \eqref{eq:L5.1}, we obtain the same conclusion \eqref{conclusion:L5.1}. 
\end{remark}

\subsection*{Asymptotics of the linear term in 2D}

Based on a refinement of the description of the asymptotic behaviour of the convolution of a certain class of functions made in \cite{Bra08}, we can control the deviation of $e^{t\Delta}(u_0,b_0)$ from $(u_0, b_0)$ in the far-field in the following sense:

\begin{lemma}
\label{lemma:6.1}
Let $f \in L^2(\R^2)$. Assume that at least one of the following conditions holds:
\begin{itemize}
	\item[(i)] Either $f(x) = o(|x|^{-3})$, or
	\item[(ii)] $\nabla f(x) = o(|x|^{-3})$, or
	\item[(iii)] $\nabla^2 f(x) = o(|x|^{-3})$.
\end{itemize}
Then, there exists a polynomial $P = P(T)$ and a function $\Tilde{\epsilon}= \Tilde{\epsilon}(x)$ such that
\begin{equation*}
	\lim_{|x| \to \infty} \Tilde{\epsilon}(x) = 0 \qquad \& \qquad \sup_{t \in (0,T)} |e^{t\Delta}f(x) - f(x)| \leq P(T)|x|^{-3}\Tilde{\epsilon}(x)
\end{equation*}
for all $T > 0$.
\end{lemma}

The proof of this lemma can be found in (\cite{Bra22}, section 6).

\subsection*{The spaces $X_{p,T}$ and $Y_{p,T}$}
Consider the functions
\begin{equation*}
	\varphi(x) = (1+|x|)\log(e+|x|)^{1/2} \quad \text{and} \quad \psi(x) = (1+|x|)^2\log(e+|x|)^{1/2}.
\end{equation*}
These weight functions are tailored to the decay hypotheses of Theorem
\ref{thm:MHDthm} and Proposition \ref{prop:MHD1} below: $\varphi$ controls the spatial rate $|x|^{-1}(\log|x|)^{-1/2}$ assumed on $u_0,b_0$, while $\psi$
controls the one-order-faster rate assumed on the gradients. For any $T > 0$ and $1 \leq p < \infty$, let
\begin{equation*}
	\|u\|_{X_{p,T}} := \esssup{t \in (0,T)}t^{1/p} \|\varphi u(t)\|{\infty} + \esssup{t \in (0,T)} t^{1/2 + 1/p} \|\psi \nabla u(t)\|_{\infty}.
\end{equation*}
We can then define the Banach space $X_{p,T}$ of measurable functions $u$ on $\R^2 \times (0,T)$ such that $\|u\|_{X_{p,T}} < \infty$. We also consider the closed subspace $Y_{p,T} \subset X_{p,T}$ defined as follows
\begin{equation*}
	\begin{split}
		Y_{p,T}=\Bigl\{u\in X_{p,T} \ | \ 
		\varphi(x)&\esssup{t \in (0,T)}  t^{1/p}|u(x,t)|\to0 \quad\text{and}\\
		& \;\psi(x)\esssup{t \in (0,T)} t^{1/p+1/2}\,|\nabla u(x,t)|
		\to0,\quad\text{as $|x|\to\infty$}\Bigr\}.
	\end{split}
\end{equation*}
We then have the following lemmas
\begin{lemma}
	\label{lem:MHD1}
	Let $B(u,v)(t) = - \int_0^t F(t-s)\ast(u \otimes v)(s) \ ds$ be the bilinear kernel operator. For all $T > 0$ and all $4 < p < \infty$, we have the estimate
	\begin{equation}
		\label{ineq:MHD1}
		\|B(u,v)\|_{X_{p,T}} \leq C_p T^{1/2 - 1/p}\left(1 + \sqrt{T} \right)\|u\|_{X_{p,T}}\|v\|_{X_{p,T}}.
	\end{equation}
	Moreover, if $u,v \in Y_{p,T_0}$, then $\BBc(u,v) \in Y_{p,T_0}$.
\end{lemma}

\begin{lemma}
	\label{lem:MHD2}
	Let $1 \leq p < \infty$ and $v \in L^p(\R^2)$, such that $v(x) = o(\varphi(x)^{-1})$ and $\nabla v(x) = o(\psi(x)^{-1})$ as $|x| \to \infty$. Then, for all $T>0$, $e^{t\Delta}v \in Y_{p,T}$ and
	\begin{equation*}
		\|e^{t\Delta}v\|_{X_{p,T}} \leq C_p(1 + T^{1 + 1/p}),
	\end{equation*}
	for some $C_p>0$ independent of $T$.
\end{lemma}
Both lemmas were proven in \cite{Bra22}.
\begin{remark}
	The spaces $X_{p,T}$ and $Y_{p,T}$ are built to track two properties
	simultaneously: the locally-in-time boundedness of the $L^{\infty}$-norm of the heat semigroup and its gradient and the spatial decay expected at spatial infinity ($\varphi^{-1}$ and $\psi^{-1}$, respectively). The closed subspace $Y_{p,T}$, in which the weighted quantities vanish as $|x|\to\infty$ rather than merely stay bounded, is what ultimately allows the little-$o$ statement of Lemma \ref{lemma:5.1} to be	invoked, and hence the sharp profile of Theorem \ref{thm:MHDthm} to be obtained.
\end{remark}

\section{Hexagonal Structure for MHD}\label{section3}

In this section we will show Theorem \ref{thm:MHDthm}. The main ingredient in such a proof is the following persistence result for MHD solutions with initial data as in the former theorem, precisely. 

\begin{proposition}
	\label{prop:MHD1}
	Let $u_0, b_0 \in L_{\sigma}^2(\R^2) \cap L^p(\R^2)$, for some $4 < p < \infty$, be such that
	\begin{equation*}
		\begin{cases}
			u_0(x) = o(|x|^{-1}(\log|x|)^{-1/2}) \\
			\nabla u_0(x) = o(|x|^{-2}(\log|x|)^{-1/2})
		\end{cases}
		\text{ and } \quad
		\begin{cases}
			b_0(x) = o(|x|^{-1}(\log|x|)^{-1/2}) \\
			\nabla b_0(x) = o(|x|^{-2}(\log|x|)^{-1/2})
		\end{cases}
	\end{equation*}
	as $|x| \to \infty$. Then, the unique Leray's solution $(u,b)$ starting from $(u_0,b_0)$ satisfies, for all $0 < T < \infty$,
	\begin{equation}
		\label{eq:MHDres1}
		\begin{split}
			&\begin{cases}
				\esssup{t \in (0,T)}t^{1/p}|u(x,t)|= o(|x|^{-1}(\log|x|)^{-1/2}) \\
				\esssup{t \in (0,T)} t^{1/2 + 1/p}|\nabla u_(x,t)|= o(|x|^{-2}(\log|x|)^{-1/2})
			\end{cases} \\
			\text{and}\hspace{30pt} \\
			&\begin{cases}
				\esssup{t \in (0,T)}t^{1/p}|b(x,t)| = o(|x|^{-1}(\log|x|)^{-1/2}) \\
				\esssup{t \in (0,T)} t^{1/2 + 1/p}|\nabla b(x,t)| = o(|x|^{-2}(\log|x|)^{-1/2})
			\end{cases}
		\end{split}
	\end{equation}
	as $|x| \to \infty$.
\end{proposition}
This means that our flow preserves the pointwise decay conditions placed on $(u_0, b_0)$, which shows that such conditions can be satisfied in physically relevant cases.

To show the existence of a solution to (\ref{eq:MHDintegral}) in such spaces (i.e., satisfying \eqref{eq:MHDres1}), we rely on the following fixed-point argument
\begin{lemma}
	\label{lem:MHDFixedPoint}
	Let $E$ be a Banach space, $\BBc$ a continuous bilinear map from $E \times E$ to $E$, and $\alpha$ a positive real number such that
	\begin{equation*}
		\alpha < \frac{1}{8\|\BBc\|}, \text{ with } \|\BBc\| := \underset{\|u\|,\|v\| \leq 1}{\sup} \|\BBc(u,v)\|.
	\end{equation*}
	For any $a, d$ in the ball $B(0,\alpha) \subset E$, there exists a unique pair $(x,y) \in B(0,2\alpha) \times B(0,2\alpha)$ satisfying
	\begin{equation}
		\label{eq:MHDFixedPoint}
		\begin{cases}
			x = a + \BBc(x,x) - \BBc(y,y) \\
			y = d + \BBc(x,y) - \BBc(y,x)
		\end{cases}
	\end{equation}
\end{lemma}
\begin{proof}
	Consider the iterative scheme
	\begin{equation*}
		\begin{cases}
			x_0 = a \quad \text{and} \quad y_0 = d, \\
			x_{n+1} = a + \BBc(x_n,x_n) - \BBc(y_n,y_n), \\
			y_{n+1} = d + \BBc(x_n,y_n) - \BBc(y_n,x_n).
		\end{cases}
	\end{equation*}
	Since $a,d \in B(0,\alpha)$, we have, via induction, that $\|x_n\|, \|y_n\| \leq 2\alpha$. Now,
	\begin{equation*}
		\begin{split}
			x_{n+1}-x_n &= [\BBc(x_n,x_n) - \BBc(x_{n-1},x_{n-1})] - [\BBc(y_n,y_n) - \BBc(y_{n-1},y_{n-1})] \\
			&= [\BBc(x_n - x_{n-1},x_n) + \BBc(x_{n-1},x_n - x_{n-1})]\\
			&- [\BBc(y_n - y_{n-1}, y_n) + \BBc(y_{n-1}, y_n - y_{n-1})],
		\end{split}
	\end{equation*}
	giving us
	\begin{equation*}
		\|x_{n+1} - x_n\| \leq 4\alpha\|\BBc\|(\|x_n - x_{n-1}\| + \|y_n - y_{n-1}\|).
	\end{equation*}
	Analogously, we have
	\begin{equation*}
		\|y_{n+1} - y_n\| \leq 4\alpha\|\BBc\|(\|x_n - x_{n-1}\| + \|y_n - y_{n-1}\|).
	\end{equation*}
	Thus, by induction on $n$,
	\begin{equation*}
		\|x_{n+1} - x_n\| + \|y_{n+1} - y_n\| \leq (8\alpha \|\BBc\|)^n(\|x_1 - x_0\| + \|y_1 - y_0\|).
	\end{equation*}
	Therefore, for all $k \geq 1$,
	\begin{equation*}
		\begin{split}
			\|x_{n+k} - x_n\| + \|y_{n+k} - y_n\| &\leq \sum_{j = 0}^{k-1} (\|x_{n+j+1} - x_{n+j}\| + \|y_{n+j+1} - y_{n+j}\|) \\
			&\leq (\|x_1 - x_0\| + \|y_1 - y_0\|) \sum_{j=0}^{k-1}(8\alpha \|\BBc\|)^{n+j} \\
			&\leq \frac{1}{1 - 8\alpha\|\BBc\|}(\|x_1 - x_0\| + \|y_1 - y_0\|)(8\alpha \|\BBc\|)^n < \epsilon,
		\end{split}
	\end{equation*}
	for $n$ large enough. Therefore, both $(x_n)_n$ and $(y_n)_n$ are Cauchy sequences. Since $E$ is a Banach space, by letting $x = \lim x_n$ and $y = \lim y_n$, and from the continuity of $\BBc$, $(x,y)$ is the desired solution. The uniqueness argument is straightforward, since, if $(x,y), (x',y')$ are both solutions to (\ref{eq:MHDFixedPoint}), then
	\begin{equation*}
		\begin{split}
			&\|x - x'\| + \|y - y'\| \leq 8 \alpha \|\BBc\| (\|x - x'\| + \|y - y'\|) \\
			\implies &(1 - 8\alpha\|\BBc\|)(\|x-x'\| + \|y-y'\|) \leq 0 \\
			\therefore\quad &x = x' \quad \text{and} \quad y = y',
		\end{split}
	\end{equation*}
	because $ 8\alpha\|\BBc\| < 1$.
\end{proof}

\begin{remark}
	Unlike the Navier-Stokes case, where the local existence in $X_{p,T}$ (or
	$Y_{p,T}$) follows from a direct Picard iteration on the single bilinear
	estimate (\ref{ineq:MHD1}), the coupled system (\ref{eq:MHDintegral}) does not reduce to a one-variable contraction because of the cross terms in each equation. Thus, any naive two-variable iteration must control
	cross terms between $u$ and $b$ at every step, and the standard
	Navier-Stokes scheme does not directly provide the needed bound on the
	sum $\|x_n\|+\|y_n\|$ that is required. Lemma \ref{lem:MHDFixedPoint}
	above isolates the abstract structure that closes this coupled
	iteration.
\end{remark}

\begin{proof}[Proof of Proposition \ref{prop:MHD1}]
	If $u_0, b_0$ satisfy the hypotheses in the proposition, then, from Lemmas \ref{lem:MHD1} and \ref{lem:MHD2}, by taking $T_0 > 0$ sufficiently small and letting $a = e^{t\Delta}u_0$ and $d = e^{t\Delta}b_0$, there exists a unique solution $(u,b) \in Y_{p,T_0} \times Y_{p,T_0}$ to (\ref{eq:MHDintegral}). Moreover, since $u_0, b_0 \in L^2_{\sigma}(\R^2)$, the iteration scheme (\ref{eq:MHDFixedPoint}) is also known to converge in $L^4\big([0,T_0], \Dot{H}^{1/2}(\R^2)\big)$ (by, possibly, making $T_0$ smaller), via a classical Fujita and Kato's argument \cite{FK64} (see also \cite{Dan05}). Therefore, our solution $(u,b)$ is in $\big(L^4_{T_0}\dot{H}^{1/2} \cap Y_{p,T_0}\big)^2$. Then, since the corresponding two-dimensional Leray solution is globally
	smooth for $t>0$ (see \cite{SM83}), strong-weak uniqueness gives that such a fixed-point solution is global, with
	\[
	\sup_{t \in [T_0,T]} (\|u(t)\|_{\infty} + \|b(t)\|_{\infty}) < \infty, \quad \forall T > T_0.
	\]
	We now work on $[T_0,T]$, for some arbitrary $T > T_0$. Let
	\begin{equation*}
		\Tilde{u}_0(x) = u(x,T_0) \qquad \text{and} \qquad \Tilde{b}_0(x) = b(x,T_0).
	\end{equation*}
	Since $u,b \in Y_{p,T_0}$, we have that 
	\begin{equation}
		\label{condT0MHD}
		\begin{cases}
			\varphi\Tilde{u}_0,~ \varphi\Tilde{b}_0 \in L^{\infty}(\R^2) \\
			\psi\nabla \Tilde{u}_0,~ \psi\nabla\Tilde{b}_0 \in L^{\infty}(\R^2)
		\end{cases}
		\qquad \text{and} \qquad
		\begin{cases}
			\varphi(x)\Tilde{u}_0(x), ~\varphi(x)\Tilde{b}_0(x) \to 0 \\
			\psi(x) \nabla\Tilde{u}_0(x), ~\psi(x)\nabla\Tilde{b}_0(x) \to 0
		\end{cases}
	\end{equation}
	as $|x| \to \infty$. For $a > 0$, let
	\begin{equation*}
		\rho_a(x) := (1 + |x|)^a
	\end{equation*}
	and set $a = 2/3$. For $T_0 \leq s \leq t \leq T$, we consider the identity
	\begin{equation}
		\label{eq:MHDfroms}
		\begin{cases}
			u(t) = e^{(t-s)\Delta}u(\cdot,s) - \int_s^t F(t-\tau) \ast ( u \otimes u - b \otimes b)(\tau) \ d\tau, \\
			b(t) = e^{(t-s)\Delta}b(\cdot,s) - \int_s^t F(t-\tau) \ast (u \otimes b - b \otimes u)(\tau) \ d\tau,
		\end{cases}
	\end{equation}
	where $F(x,t)$ is the kernel of the operator $e^{t\Delta}\P\dive$. Now, from \cite{Vig05}, we know that, if $v = o(\varphi(x)^{-1})$, as $|x| \to \infty$, then
	\begin{equation*}
		\sup_{t \in [0, T]} \|\rho_ae^{t\Delta}v\|_{\infty} \leq A\|\rho_av\|_{\infty} \qquad \text{and} \qquad \sup_{t \in [0, T]} \|\rho_a \nabla e^{t\Delta}v\|_{\infty} \leq A\|\rho_a\nabla v\|_{\infty},
	\end{equation*}
	for some constant $A = A(T)$. Combining these estimates with equation (\ref{eq:MHDfroms}), we get, for all $T_0 \leq s \leq t \leq T$,
	\begin{equation*}
		\begin{split}
			\|\rho_a u(t)&\|_{\infty} + \|\rho_a \nabla u(t)\|_{\infty} \\
			&\leq A(\|\rho_a u(s)\|_{\infty} + \|\rho_a \nabla u(s)\|_{\infty}) \\
			&+ 2B(t-s)^{1/2}\left( \sup_{\tau \in [s,t]} \|u(\tau)\|_{\infty} \right)\left( \sup_{\tau \in [s,t]} \|\rho_au(\tau)\|_{\infty} +\sup_{\tau \in [s,t]} \|\rho_a \nabla u(\tau)\|_{\infty} \right) \\
			&+ 2B(t-s)^{1/2}\left( \sup_{\tau \in [s,t]} \|b(\tau)\|_{\infty} \right)\left( \sup_{\tau \in [s,t]} \|\rho_a b(\tau)\|_{\infty} +\sup_{\tau \in [s,t]} \|\rho_a \nabla b(\tau)\|_{\infty} \right).
		\end{split}
	\end{equation*}
	and
	\begin{equation*}
		\begin{split}
			\|\rho_a b(t)&\|_{\infty} + \|\rho_a \nabla b(t)\|_{\infty} \\
			&\leq A(\|\rho_a b(s)\|_{\infty} + \|\rho_a \nabla b(s)\|_{\infty}) \\
			&+ 2B(t-s)^{1/2}\left( \sup_{\tau \in [s,t]} \|u(\tau)\|_{\infty} \right)\left( \sup_{\tau \in [s,t]} \|\rho_ab(\tau)\|_{\infty} +\sup_{\tau \in [s,t]} \|\rho_a \nabla b(\tau)\|_{\infty} \right) \\
			&+ 2B(t-s)^{1/2}\left( \sup_{\tau \in [s,t]} \|b(\tau)\|_{\infty} \right)\left( \sup_{\tau \in [s,t]} \|\rho_a u(\tau)\|_{\infty} +\sup_{\tau \in [s,t]} \|\rho_a \nabla u(\tau)\|_{\infty} \right).
		\end{split}
	\end{equation*}
	Now, let
	\begin{equation*}
		\begin{split}
			&\alpha(v;t) := \|\rho_a v(t)\|_{\infty} + \|\rho_a \nabla v(t)\|_{\infty}, \quad
			\beta(s,t) := \sup_{\tau \in [s,t]}\|u(\tau)\|_{\infty} + \sup_{\tau \in [s,t]} \|b(\tau)\|_{\infty},\\
			&\gamma(v;s,t) := \sup_{\tau \in [s,t]}\|\rho_a v(\tau)\|_{\infty} + \sup_{\tau \in [s,t]} \|\rho_a\nabla v(\tau)\|_{\infty}.
		\end{split}
	\end{equation*}
	Combining both previous expressions, we get
	\begin{equation*}
		\begin{split}
			\alpha(u;t) + \alpha(b;t) \leq A[\alpha(u;s) + \alpha(b;s)] + 2B(t-s)^{1/2}\beta(s,t)[\gamma(u;s,t) + \gamma(b;s,t)].
		\end{split}
	\end{equation*}
	Now, let's assume $(u,b) \not \equiv (0,0)$ on $[T_0,T]$. Starting from $T_0$, consider the strictly increasing sequence of times $(T_n)_{n \geq 0}$ such that, for $n \geq 0$,
	\begin{equation*}
		2B (T_{n+1} - T_n)^{1/2} \sup_{\tau \in [T_0,T]} \Big(\|u(\tau)\|_{\infty} + \|b(\tau)\|_{\infty}\Big)= 1/2.
	\end{equation*}
	Let $N \in \N$ be such that $T_N \leq T < T_{N+1}$. We have
	\begin{equation*}
		N \leq (T-T_0)\left[4B \sup_{\tau \in [T_0,T]}\Big(\|u(\tau)\|_{\infty} + \|b(\tau)\|_{\infty}\Big)\right]^2 < N+1.
	\end{equation*}
	For $n \leq N$, consider the interval $\Delta_n = [T_0,T] \cap [T_n, T_{n+1}]$, and set
	\begin{equation*}
		M_n = \sup_{\tau \in \Delta_n} \Big(\|\rho_a u(\tau)\|_{\infty} + \|\rho_a b(\tau)\|_{\infty}\Big) + \sup_{\tau \in \Delta_n} \Big(\|\rho_a \nabla u(\tau)\|_{\infty} + \|\rho_a \nabla b(\tau)\|_{\infty}\Big).
	\end{equation*}
	From (\cite{Vig05}, Corollary 5) and our previous estimates,
	notice that, for each $n \leq N$ and $T_n \leq t \leq T_{n+1}$,
	\begin{equation*}
		\sup_{\tau \in \Delta_n} \Big(\|\rho_a u(\tau)\|_{\infty} + \|\rho_a b(\tau)\|_{\infty}\Big) \leq A(\|\rho_a u(T_n)\|_{\infty} + \|\rho_a b(T_n)\|_{\infty})
		+ \frac{1}{2} \sup_{\tau \in \Delta_n} \Big(\|\rho_a u(\tau)\|_{\infty} + \|\rho_a b(\tau)\|_{\infty}\Big)
	\end{equation*}
	and
	\begin{equation*}
		\begin{split}
			\sup_{t \in \Delta_n}\Big(\|\rho_a \nabla u(t)\|_{\infty} + \|\rho_a \nabla b(t)\|_{\infty}\Big) &\leq A\Big(\|\rho_a \nabla u(T_n)\|_{\infty} + \|\rho_a \nabla b(T_n)\|_{\infty} \Big) \\
			&+ \frac{1}{2} \sup_{\tau \in \Delta_n} \Big(\|\rho_a \nabla u(\tau)\|_{\infty} + \|\rho_a \nabla b(\tau)\|_{\infty} \Big)
		\end{split}
	\end{equation*}
	which implies
	\begin{equation*}
		M_n \leq 2A M_{n-1}, \quad 1 \leq n \leq N.
	\end{equation*}
	Analogously, for $n = 0$, we get
	\begin{equation*}
		M_0 \leq 2A\Big( \|\rho_a \Tilde{u}_0\|_{\infty} + \|\rho_a \Tilde{b}_0\|_{\infty} + \|\rho_a \nabla\Tilde{u}_0\|_{\infty} + \|\rho_a \nabla\Tilde{b}_0\|_{\infty} \Big).
	\end{equation*}
	Therefore, taking the product over $n = 0, \ldots, N$,
	\begin{equation*}
		\begin{split}
			\sup_{\tau \in [T_0,T]} \Big( \|\rho_a u(\tau)\|_{\infty} + \|\rho_a b(\tau)\|_{\infty} \Big) &+ \sup_{\tau \in [T_0,T]} \Big( \|\rho_a \nabla u(\tau)\|_{\infty} + \|\rho_a \nabla b(\tau)\|_{\infty} \Big) \\
			&\leq (2A)^{N+1}\Big( \|\rho_a \Tilde{u}_0\|_{\infty} + \|\rho_a \Tilde{b}_0\|_{\infty} + \|\rho_a \nabla\Tilde{u}_0\|_{\infty} + \|\rho_a \nabla\Tilde{b}_0\|_{\infty} \Big)
		\end{split}
	\end{equation*}
	But, from our choice of $N$,
	\begin{equation*}
		(2A)^{N+1} = 2A \exp(N\log(2A)) \leq A_1 \exp\left(A_2 \sup_{t \in [T_0,T]} \Big(\|u(t)\|_{\infty} + \|b(t)\|_{\infty}\Big)^2\right),
	\end{equation*}
	where $A_1$ and $A_2$ depend only on $T$, and are locally bounded functions of $T$. We then conclude that
	\begin{equation*}
		\begin{split}
			&\sup_{t \in [T_0,T]} \Big( \|\rho_a u(t)\|_{\infty} + \|\rho_a b(t)\|_{\infty} \Big) + \sup_{t \in [T_0,T]} \Big( \|\rho_a \nabla u(t)\|_{\infty} + \|\rho_a \nabla b(t)\|_{\infty} \Big) \\
			&\leq A_1 \exp\Big(A_2 \sup_{t \in [T_0,T]} (\|u(t)\|_{\infty} + \|b(t)\|_{\infty})^2\Big) \Big(\|\rho_a \Tilde{u}_0\|_{\infty} + \|\rho_a \Tilde{b}_0\|_{\infty} + \|\rho_a \nabla\Tilde{u}_0\|_{\infty} + \|\rho_a \nabla\Tilde{b}_0\|_{\infty} \Big).
		\end{split}
	\end{equation*}
	Thus, from the above estimates,
	\begin{equation*}
		\sup_{t \in [T_0,T]} \Big( |u(x,t)| + |b(x,t)| \Big) + \sup_{t \in [T_0,T]} \Big( |\nabla u(x,t)| + |\nabla b(x,t)| \Big) \leq C(\rho_{2/3})^{-1},
	\end{equation*}
	where $C$ depends only on $T_0$, $T$ and $(u_0,b_0)$. Now, from (\ref{eq:MHDfroms}), with $s = T_0$,
	\begin{equation*}
		\begin{cases}
			u(x,t) = e^{(t-T_0)\Delta}\Tilde{u}_0(x) - \int_{T_0}^t F(t-\tau) \ast ( u \otimes u - b \otimes b)(\tau) \ d\tau \\
			b(x,t) = e^{(t-T_0)\Delta}\Tilde{b}_0(x) - \int_{T_0}^t F(t-\tau) \ast (u \otimes b - b \otimes u)(x,\tau) \ d\tau
		\end{cases}\quad.
	\end{equation*}
	For the quadratic term, 
	\begin{equation*}
		\begin{split}
			\sup_{t \in [T_0, T]} \int_{T_0}^t |F(t-s)|\ast(|u \otimes u| + |b \otimes b|)(x,s) \ \d s
			\leq C(T-T_0)^{1/2}(\rho_{4/3}(x))^{-1} 
			= o(\varphi(x)^{-1}),
		\end{split}
	\end{equation*}
	as $|x| \to \infty$. As for the linear term, recalling (\ref{condT0MHD}), we have $\varphi\, e^{(t-T_0)\Delta}\Tilde{u}_0 \in L^{\infty}(\R^2)$ and $\sup_{t \in [T_0,T]} |e^{(t-T_0)\Delta}\Tilde{u}_0|(x) = o(\varphi(x)^{-1})$, as $|x| \to \infty$. Hence,
	\begin{equation*}
		\varphi \, \sup_{t \in [T_0, T]} |u(\cdot,t)| \in L^{\infty}(\R^2) \qquad \text{and} \qquad \varphi(x) \, \sup_{t \in [T_0,T]} |u(x,t)| \to 0, \quad \text{as } |x| \to \infty.
	\end{equation*}
	Analogously, 
	\begin{equation*}
		\varphi \, \sup_{t \in [T_0, T]} |b(\cdot,t)| \in L^{\infty}(\R^2) \qquad \text{and} \qquad \varphi(x) \, \sup_{t \in [T_0,T]} |b(x,t)| \to 0, \quad \text{as } |x| \to \infty.
	\end{equation*}
	For the gradient estimates, after three iterations of a similar bootstrapping procedure, we get
	\begin{equation*}
		\begin{cases}
			\psi \, \sup_{t \in [T_0, T]} |\nabla u(\cdot,t)| \in L^{\infty}(\R^2) \qquad \text{and} \qquad \psi(x) \, \sup_{t \in [T_0,T]} |\nabla u(x,t)| \to 0 \\
			\psi \, \sup_{t \in [T_0, T]} |\nabla b(\cdot,t)| \in L^{\infty}(\R^2) \qquad \text{and} \qquad \psi(x) \, \sup_{t \in [T_0,T]} |\nabla b(x,t)| \to 0
		\end{cases}\, ,
	\end{equation*}
	as $|x| \to \infty$. This concludes the proof of our proposition.
\end{proof}
Now, we can focus on Theorem \ref{thm:MHDthm}.
\begin{proof}[Proof of Theorem \ref{thm:MHDthm}]
	First, let's prove the case for the velocity field. From our assumptions, $u$ and $b$ satisfy (\ref{eq:MHDres1}). Therefore, for $\kappa = 1/2 + 2 /p < 1$, both $u \otimes u$ and $b \otimes b$ satisfy (\ref{eq:L5.1}). Thus, we can write
	\begin{equation*}
		\begin{split}
			u(x,t) &= e^{t\Delta}u_0 - \LL(u \otimes u - b \otimes b)(x,t)  \\
			&= e^{t\Delta}u_0(x) -\mathfrak{F}(x) : \int_0^t \int (u \otimes u - b \otimes b)(y,s) \dy \d s + |x|^{-3}\epsilon(x,t) \\
			&= e^{t\Delta}u_0(x) -\nabla H(x,t) + o_t(|x|^{-3}),
		\end{split}
	\end{equation*}
	where, following the decomposition (\ref{def:KernelDecomPBrandolese}),
	\begin{equation*}
		\label{def:MHD-HDefinition}
		\begin{split}
			H(x,t)
			&:=\Bigl[\nabla^2 E_2(x)\colon\int_0^t\!\!\int (u\otimes u - b \otimes b)(y,s) \ dy \hspace{1pt} ds\Bigr]\\
			&=\sum_{h,k=1}^2\partial_h\partial_k E_2(x)\int_0^t\!\!\int (u_hu_k - b_hb_k)(y,s) \ dy \hspace{1pt} ds.
		\end{split}
	\end{equation*}
	If we now let
	\begin{equation*}
		\begin{pmatrix}
			{a} & {c} \\
			{c} & {d}
		\end{pmatrix}
		:=
		\begin{pmatrix}
			\int_0^t\int (u_1^2 - b_1^2)(y,s) \ dy \hspace{1pt} ds & \int_0^t\int 2(u_1u_2 - b_1b_2)(y,s) \ dy \hspace{1pt} ds \\ 
			\int_0^t\int 2(u_1u_2 - b_1b_2)(y,s) \ dy \hspace{1pt} ds & \int_0^t\int (u_2^2 - b_2^2)(y,s) \ dy \hspace{1pt} ds
		\end{pmatrix},
	\end{equation*}
	a direct calculation gives us
	\begin{equation*}
		\begin{split}
			\nabla {H}(R,\theta,t) &= \frac{\sqrt{({d} - {a})^2+{c}^2}}{\pi R^3}\begin{pmatrix}
				\cos(3\theta + \phi(t)) \\
				\sin(3\theta + \phi(t))
			\end{pmatrix} \\
			&= \frac{{L}(t)}{R^3}\begin{pmatrix}
				\cos(3\theta + \phi(t)) \\
				\sin(3\theta + \phi(t))
			\end{pmatrix}
		\end{split}
	\end{equation*}
	where $x = R(\cos(\theta), \sin(\theta))$ and, for each $t > 0$, the angle $\phi(t) \in [0,2\pi)$ is uniquely defined by the system
	\begin{equation*}
		\begin{cases}
			\cos \phi(t) = \frac{{d} - {a}}{\sqrt{({d} - {a})^2+{c}^2}}\\
			\sin \phi(t) = \frac{{c}}{\sqrt{({d} - {a})^2+{c}^2}}
		\end{cases} \, ,
	\end{equation*}
	for all $t > 0$ such that $L(t) \neq 0$ (for the cases where $L(t) = 0$, one may simply set $\phi(t) = 0$).
	As for the magnetic field, the same argument lets us write, as $|x| \to \infty$,
	\begin{equation*}
		\begin{split}
			b(x,t) &= e^{t\Delta}b_0(x) - \mathfrak{F}(x) : \int_0^t\int (u \otimes b - b \otimes u)(y,s) \dy \d s + o_t(|x|^{-3}) \\
			&= e^{t\Delta}b_0(x) - \nabla G(x,t) + o_t(|x|^{-3}),
		\end{split}
	\end{equation*}
	where, in the same spirit of the first case, $G(x,t)$ is given by
	\begin{equation*}
		\begin{split}
			G(x,t) &:= \left[ \nabla^2E_2(x) : \int_0^t \int (u \otimes b - b \otimes u)(y,s) \dy \d s \right] \\
			&= \sum_{h,k=1}^2 \partial_h\partial_k E_2(x) \int_0^t \int (u_hb_k - u_kb_h)(y,s) \dy \d s \\
			&= 0,
		\end{split}
	\end{equation*}
	since $\partial_h \partial_k E_2(x) = \partial_k\partial_h E_2(x)$. Therefore,
	\begin{equation*}
		b(x,t) = b_0(x) + o_t(|x|^{-3}), \qquad \text{as } |x| \to \infty,
	\end{equation*}
	as desired.
\end{proof}
\begin{remark}
	Notice that the condition $\kappa = \frac{1}{2} + \frac{2}{p} < 1$, needed for Lemma \ref{lemma:5.1}, is equivalent to $p > 4$. Hence, the extra $L^p$ hypothesis placed on the initial datum $(u_0, b_0)$ in Theorem \ref{thm:MHDthm} and Proposition \ref{prop:MHD1} is not arbitrary.
\end{remark}

	\appendix
	\section{Where do they begin?}
	A natural question that may arise from this result is: how far is ``far enough" to get the hexagonal structures? By considering a slight improvement on the spatial decay of the initial data, it is possible to estimate that.
	Suppose that, n addition to the other conditions in \eqref{eq:InitialDecay}, we assume that
	\begin{equation*}
		|u_0(x)|, |b_0(x)| = o(|x|^{-1-\alpha}), \quad \text{as} \quad |x| \to \infty,
	\end{equation*}
	for some $\alpha \in (\frac{1}{2},2]$. When $\alpha =2$ we are in the {\it spatially well-localized case} defined in \cite{Bra22}. Now, by a variation of Lemma \ref{lemma:5.1} (see remark \ref{needforappendix}), we can find an upper bound for $\epsilon(x,t)$ in \eqref{conclusion:L5.1}. Namely,
	\begin{equation*}
		|\epsilon(x,t)| \leq Ct^{1/2 - 1/p}\delta(x), \qquad \delta(x) \to 0 \text{ as } |x| \to \infty.
	\end{equation*}
	
	To ensure that $|\epsilon(x,t)|$ is relatively small with respect to the $|x|^{-3}$ term, we define, for $t > 0$ such that $L(t) > 0$,
	\begin{equation*}
		R(t) := \inf \left\{R>0 \, : \, \delta(x) \leq \frac{1}{4}t^{1/p-1/2}L(t), \quad \forall|x| \geq R \right\}.
	\end{equation*}
	If we then consider initial data $u_0,b_0\in L^2_\sigma (\R^2)$ such that the global solution $(u,b)$ satisfies  $L(t) \geq c > 0$, for $t \geq t_0$, for some $t_0>0$, then, from \eqref{eq:uProfile},
	\begin{equation}\label{difference}
		\|u(x,t) - e^{t\Delta}u_0\|_2^2 \geq C R(t)^{-4}.
	\end{equation}
	For the Navier-Stokes case, we can use Wiegner's theorem \cite{Wie87} to estimate from above the difference \eqref{difference}. The same (Wiegner's result) can be done for the MHD system, see e.g \cite{wiegnerMHDLiu}, by estimating the decay of $\|e^{t\Delta}u_0\|_2$ and $\|e^{t\Delta}b_0\|_2$. Doing so, we obtain that, for $t \geq t_0$,
	\begin{itemize}
		\item[(i)] If $\alpha \in (\frac{1}{2},1]$, then
		\[
		R(t) \geq C(1+t)^{\beta/2}, \qquad \forall\, 0 < \beta<\alpha.
		\]
		\item[(ii)] If $\alpha \in (1, 2]$, then
		\[
		R(t) \geq C\sqrt{1+t}.
		\]
	\end{itemize}
	\vspace{.5cm}

\scriptsize{
	\textbf{Acknowledgments.}  C. F. Perusato was partially supported by CNPq grant No. 305021/2026-5. M. Nasi was partially supported by CAPES-PDSE fellowship No. 8881.219946/2025-01 and CAPES No. 88887.149045/2025-00. The authors would like to express their gratitude to L. Brandolese and P. Zingano for their helpful suggestions and fruitful discussions concerning this work.

\end{document}